\documentclass [12pt,oneside,a4paper,mathscr]{amsart}

\usepackage {amscd,amsmath,amssymb,euscript}

\newtheorem{thm}{Theorem}[section]
\newtheorem{cor}[thm]{Corollary}
\newtheorem{lemma}[thm]{Lemma}
\newtheorem{prop}[thm]{Proposition}

\theoremstyle{definition}

\theoremstyle{remark}
\newtheorem*{pf}{Proof}

\newtheorem{rmk}[thm]{Remark}

\renewcommand{\qed}{\ifhmode\unskip\nobreak\fi\quad\ensuremath\square}

\newcommand{\col}[2]{\left( \begin{array}{c} #1 \\ #2 \end{array} \right)}
\newcommand{\tensor}{\otimes}
\newcommand{\M}{\operatorname{\mathcal
M}}
\newcommand{\ha}{\frac{1}{2}}
\newcommand{\ts}{2 *}
\newcommand{\lb}{\langle}
\newcommand{\rb}{\rangle}

\newcommand{\Jac}{\operatorname{Jac}}
\newcommand{\Se}{\operatorname{S}}
\newcommand{\bua}{\bigg{\uparrow}}
\newcommand{\embed}{\hookrightarrow}
\renewcommand{\P}{\mathcal{P}}
\newcommand{\dual}{\vee}
\newcommand{\mat}[4]{\left( \begin{array}{cc} #1 & #2 \\ #3 & #4
\end{array} \right)}
\newcommand{\K}{K}

\newcommand{\m}{{\lambda}}
\renewcommand{\H}{\operatorname{H}}
\newcommand{\J}{\operatorname{J}}
\newcommand{\eu}{\operatorname{\chi}}
\newcommand{\Phit}{\widetilde{\Phi}}
\newcommand{\Psit}{\widetilde{\Psi}}
\newcommand{\dg}{\operatorname{d}}
\newcommand{\rk}{\operatorname{r}}

\newcommand{\xt}{\tilde{x}}
\newcommand{\yt}{\tilde{y}}

\newcommand{\Xt}{\widetilde{X}}
\newcommand{\Yt}{\widetilde{Y}}
\newcommand{\C}{\mathbb C}
\newcommand{\cl}{\operatorname{c}}
\newcommand{\ch}{\operatorname{ch}}
\newcommand{\Q}{\operatorname{\mathbb Q}}

\newcommand{\NS}{\operatorname{NS}}
\newcommand{\T}{\operatorname{T}}

\newcommand{\Hom}{\operatorname{Hom}}
\newcommand{\lHom}{\operatorname{{\mathcal H}om}}
\newcommand{\lra}{\longrightarrow}
\newcommand{\lRa}[1]{\stackrel{#1}{\lra}}
\newcommand{\lla}{\longleftarrow}
\newcommand{\lLa}[1]{\stackrel{\,#1}{\lla}}
\newcommand{\bda}{\bigg{\downarrow}}
\newcommand{\isom}{\cong}
\newcommand{\Z}{{\mathbb Z}}
\newcommand{\PP}{{\mathbb P}}
\newcommand{\Ft}{\tilde{F}}
\newcommand{\QQ}{{\mathcal Q}}
\newcommand{\EE}{{\mathcal E}}
\newcommand{\FF}{{\mathcal F}}

\newcommand{\D}{\operatorname{D}}
\newcommand{\R}{\mathbf R}

\newcommand{\OO}{{\mathscr O}}
\newcommand{\Ltensor}{\mathbin{\buildrel{\mathbf L}\over{\tensor}}}
\newcommand{\Ext}{\operatorname{Ext}}

\newcommand{\Sym}{\operatorname{Sym}}
\newcommand{\Num}{\operatorname{Num}}

\begin{document}
\normalsize
\title[]{Complex surfaces with equivalent derived categories}
\author[]{Tom Bridgeland \and Antony Maciocia}

\begin{abstract}
We examine the extent to which a smooth minimal
complex projective
surface $X$ is determined by its
derived category of coherent sheaves $\D(X)$. To do this we find, for each such surface
$X$, the set of surfaces $Y$ for which there
exists a Fourier-Mukai transform $\D(Y)\to \D(X)$.
\end{abstract}

\maketitle

\section{Introduction}

This paper addresses the question: to what extent is a smooth projective variety $X$ determined by its bounded derived category
of coherent
sheaves $\D(X)$? Recall \cite{ha1} that $\D(X)$ is a
triangulated category, whose objects are bounded complexes of coherent
sheaves on $X$. If $Y$ is another smooth projective variety, an
equivalence of categories
$$\Phi\colon\D(Y)\lra\D(X)$$
preserving the triangles is called a Fourier-Mukai
transform.
Put another way then, our
problem is to
find, for a given variety $X$, the set of \emph{Fourier-Mukai partners}
of $X$, i.e. the set of varieties $Y$ for which there
exists a Fourier-Mukai transform relating $X$ and $Y$.

\smallskip

This problem is interesting for several reasons. Firstly,
Fourier-Mukai (FM) transforms have shown themselves to be important tools for  studying
moduli spaces of
sheaves \cite{br1}, \cite{mac}, \cite{muk4}, and it is therefore natural to attempt to
classify them. Secondly, the theory of Fourier-Mukai-type transforms
seems to provide the correct language for describing certain
geometrical dualities
suggested by string theory. As a particular example of this,  M. Kontsevich's
homological mirror conjecture \cite{kon} predicts that all mirror varieties of a given
variety have equivalent
derived categories. Thus, the existence of distinct FM partners of a
variety $X$ may relate to the possibility that the conjectural mirror map is not a
well-defined bijection at $X$.

\smallskip

The first example of a non-trivial FM transform was given by S. Mukai
in 1981 and related the
derived category of an abelian variety with the derived category of
the dual variety \cite{muk2}. Since then further examples have been given, involving K3
surfaces \cite{bbhr}, \cite{br2}, abelian surfaces \cite{mac},
elliptic surfaces \cite{br1} and Enriques and bielliptic surfaces
\cite{br3}. Clearly, some sort of classification is
in order.

\smallskip

The classification of FM transforms splits naturally into two
parts. Given a smooth projective variety $X$ these are

\begin{enumerate}
\item[(a)]
find the set of FM partners of
$X$, that is, the set of varieties $Y$ for which there exists a FM
transform $\D(Y)\to\D(X)$,
%\smallskip
\item[(b)]
find the
group of FM transforms $\D(X)\to\D(X)$.
\end{enumerate}
%We briefly review what is known about each of these problems.

\smallskip

When $X$ has
ample canonical or anticanonical bundle a complete solution was obtained
by
A. Bondal and D. Orlov \cite{bo1}, \cite{bo2}. In this case the answer is rather trivial, in
that the only FM partner of $X$ is $X$ itself, and all
autoequivalences of $\D(X)$ are generated by shifts, automorphisms of
$X$ and twists by line bundles.

Remarkably, Orlov also managed to solve both problems when $X$ is an
abelian variety. In this case the solution is very interesting and
highly non-trivial
\cite{orl2}.

These two results together give a simple answer to both problems in the
case when $X$ has
dimension one. In particular, it is possible to prove that the only FM
partner of
a curve $X$ is $X$ itself.

%\smallskip

In this paper
we solve Problem (a) for minimal complex surfaces. We obtain the following theorem, which
will be explained in greater detail below.

\begin{thm}
\label{chin}
Let $X$ be a smooth minimal complex projective surface, and let $Y$ be a Fourier-Mukai
partner of $X$, not isomorphic to $X$.  Then either $X$ is an
elliptic surface, and $Y$ is another elliptic
surface obtained as in \cite{br1}, \cite{fr}, by taking a relative
Picard scheme of the elliptic fibration on $X$, or $X$ is of K3 or abelian type, and $Y$
is of the same type, with Hodge-isometric transcendental lattice.
\end{thm}

\begin{cor}
The number of FM partners of a smooth minimal complex projective
surface is finite.
\end{cor}

The proof of Theorem \ref{chin} is rather long, since each different type of
surface appearing in the Enriques classification must be analysed
separately. For surfaces of Kodaira dimension 0, the problem is mostly
lattice-theoretic, and we rely heavily on results of V. Nikulin. Other
surfaces are best treated with more geometric methods. In particular,
it becomes important to classify
curves with non-positive self-intersection which do not intersect the
canonical divisor.

Problem (b) for surfaces is much more difficult. In particular,
determining the group of autoequivalences of the derived category of a K3 surface
seems to be of considerable interest.

\subsection*{Notation}
All varieties will be over $\C$. Given a variety $X$, the translation (or shift) functor on $\D(X)$ is written $[1]$, so that the
symbol $E[m]$ means the object $E$ of $\D(X)$ shifted to the left by
$m$ places. By a sheaf on $X$ we mean a coherent $\OO_X$-module, and a point of $X$ always means a closed (or
geometric) point. The structure sheaf of such a point $x\in X$ will be
denoted $\OO_x$. The canonical bundle of a smooth projective variety
$X$ is denoted $\omega_X$. By a lattice we mean a free abelian group
of finite rank with a non-degenerate $\Z$-valued symmetric bilinear form.

% ************************************************************************

\section{Preliminaries on Fourier-Mukai transforms}

Throughout this section we fix a pair of smooth projective
varieties $X$ and $Y$.

\subsection{}
A Fourier-Mukai
transform relating $X$ and $Y$ is an exact\footnote{A functor between
triangulated categories is \emph{exact} if it commutes with the
translation functors, and takes distinguished triangles to
distinguished triangles.}
equivalence of categories
$$\Phi\colon\D(Y)\lra \D(X).$$
Due to a theorem of
Orlov \cite{orl1}, it is known that for any such equivalence $\Phi$
there is an object
$\P$ of $\D(Y\times X)$, unique up to isomorphism, such that $\Phi$ is
isomorphic to the functor defined by the formula
$$\Phi^{\P}_{Y\to X}(-)=\R\pi_{X,*}(\P\Ltensor\pi_Y^*(-)),$$
where $Y\lLa{\,\pi_Y}Y\times X\lRa{\pi_X} X$ are the projection maps.
The object $\P$ is called the
\emph{kernel} of the transform $\Phi$.

We say that $X$ and $Y$ are \emph{Fourier-Mukai partners} if there is
a FM transform relating $X$ and $Y$. This is equivalent to the
statement that $\D(X)$ and $\D(Y)$ are equivalent as triangulated
categories.

\begin{lemma}
\label{extra}
If $X$ and $Y$ are FM partners then $\dim(X)=\dim(Y)$ and the
canonical line bundles $\omega_X$ and $\omega_Y$ have the same order.
\end{lemma}

\begin{pf}
Define the \emph{Serre functor} $\Se_X$ on the category $\D(X)$ by the formula
$$\Se_X(-)=\omega_X\tensor(-)[\dim X].$$
In \cite[Prop. 3.4]{boka} it is shown that any FM transform
$$\Phi\colon\D(Y)\lra\D(X)$$
commutes
with  the Serre functors on $X$ and $Y$. Thus if $\Psi$ is a quasi-inverse to
the equivalence
$\Phi$,
there is an
isomorphism of functors
$$\Se_Y\isom \Psi \circ \Se_X\circ \operatorname{\Phi}.$$
The lemma is an immediate consequence of this.
\qed
\end{pf}

Given a FM transform $\Phi\colon\D(Y)\to \D(X)$, and an object $E$ of
$\D(Y)$, let us write
$$\Phi^i(E)=\H^i(\Phi(E))$$
for the $i$th cohomology sheaf of the object $\Phi(E)$ of $\D(X)$. We shall call $\Phi$ a \emph{sheaf transform} if there
is an integer $p$ such that for
each point $y\in Y$,
$$\Phi^i(\OO_y)=0\mbox{ unless }i=p.$$
An equivalent condition, \cite[Lemma 4.3]{br2}, is that the kernel of $\Phi$ is
concentrated in some degree $p$, and is flat over $Y$.

\subsection{}
Let $E$ and $F$ be objects of $\D(Y)$. For each integer $i$ one
defines a vector space
$$\Hom^i_{\D(Y)}(E,F)=\Hom_{\D(Y)}(E,F[i]).$$
Recall that if $E$ and $F$ are concentrated in degree $0$ then these spaces
are just the $\Ext$-groups of the sheaves $E$ and $F$, i.e.
$$\Hom^i_{\D(Y)}(E,F)=\Ext^i_Y(E,F).$$
The following trivial but useful observation is sometimes referred to
as the \emph{Parseval theorem}.

\begin{lemma}
\label{percy}
Let $\Phi\colon\D(Y)\to\D(X)$ be a FM transform, and take objects $E$
and $F$ of $\D(Y)$.  Then
$$\Hom^i_{\D(X)}(\Phi(E),\Phi(F))=\Hom^i_{\D(Y)}(E,F).$$
\end{lemma}

\begin{pf}
Immediate because $\Phi$ is an equivalence of categories, commuting with the
translation functors in $\D(Y)$ and $\D(X)$.
\qed
\end{pf}

The lemma implies that
$$\eu(\Phi(E),\Phi(F))=\eu(E,F),$$
where $\eu(E,F)$ denotes the relative Euler character
$$\eu(E,F)=\sum_i (-1)^i\dim \Hom^i_{\D(Y)}(E,F).$$
This relative Euler character is given in terms of the Chern characters\footnote{The Chern character of an object of the derived
category is just  the alternating sum of the Chern characters of the
cohomology sheaves.}
of $E$ and $F$ by the Riemann-Roch
theorem. For example, if $Y$ is a surface, then
\begin{align*}
\eu(E,F)=&\rk(E)\ch_2(F)-\cl_1(E)\cdot
\cl_1(F)+\rk(F)\ch_2(E) \\
&+\frac{1}{2}(\rk(F)\cl_1(E)-\rk(E)\cl_1(F))\cdot
\K_Y+\rk(E)\rk(F)\eu(\OO_Y),
\end{align*}
where $\K_Y$ is the
first Chern class of the canonical line bundle $\omega_Y$. In
particular, if $E$ and $F$ are torsion sheaves
$$\eu(E,F)=-\cl_1(E)\cdot\cl_1(F).$$

\subsection{}
\label{cc}
Grothendieck's
Riemann-Roch theorem implies that for any FM transform $\Phi\colon\D(Y)\to \D(X)$
there is a linear map of vector spaces
$$\phi\colon \H^*(Y,\Q)\lra \H^*(X,\Q)$$
making the following diagram commute
$$
\begin{array}{ccc}
\D(Y) &\lRa{\Phi} &\D(X) \\
\scriptstyle{\ch}\bda  &&\scriptstyle{\ch}\bda \\
\H^*(Y,\Q) &\lRa{\phi} & \H^*(X,\Q),
\end{array}
$$
where $\ch$ denotes the operation of taking the Chern character.

The
proof of \cite[Theorem 4.9]{muk5} shows that $\phi$ is an isomorphism
of vector spaces. Furthermore,  since $\phi$ is given by an algebraic class on the product $Y\times X$,
it preserves the parity of the degree of cohomology
classes. One therefore has

\begin{prop}
\label{picard}
Surfaces with equivalent derived categories have the same Picard
number, and the same topological Euler number.
\qed
\end{prop}

\subsection{}
An important property of FM transforms is that they preserve families
of sheaves. Let $\Phi\colon\D(Y)\to
\D(X)$ be a FM transform, take a scheme $S$ of finite type over $\C$, and let $\EE$ be a sheaf on $S\times Y$, flat
over $S$.

\begin{prop}
\label{fam}
The set of points $s\in S$ for which
the object $\Phi(\EE_s)$ of $\D(X)$ is concentrated in degree 0 is the
set of points of an open subset $U$ of $S$ (possibly
empty). Furthermore there is a sheaf $\FF$ on $U\times X$, flat over
$U$, such that for any point $s\in U$, $\FF_s=\Phi(\EE_s)$.
\end{prop}

\begin{pf}
See \cite[Chapter 6]{thesis} or \cite[Theorem 1.6]{muk4}.
\qed
\end{pf}

As a consequence one has

\begin{lemma}
\label{birat}
Let $\Phi\colon\D(Y)\to \D(X)$ be a FM
transform, and suppose there is a point $y\in Y$, such that
$$\Phi(\OO_y)=\OO_x[p],$$
for some point $x\in X$ and some integer $p$. Then $X$ and $Y$ are
birationally equivalent.
\end{lemma}

\begin{pf}
By Prop. \ref{fam} there is an open subset $V\subset Y$ such that for each point
$y\in V$, there is a point $f(y)\in X$ with
$$\Phi(\OO_y)=\OO_{f(y)}[p].$$

The kernel of $\Phi$, restricted to $V\times X$, is supported on the graph of $f$,
so $f$ is a morphism $V\to X$, and hence defines a
birational map $Y\dashrightarrow X$. Since $\Phi$ is an
equivalence this birational map has an inverse, so $X$ and $Y$ are
birationally equivalent.
\qed
\end{pf}

\begin{rmk}
\label{ratty}
Suppose the conditions of Lemma \ref{birat} hold, and that $X$
is a minimal surface of non-negative Kodaira dimension. Then $Y$ must be the blow-up of $X$ at $r\geq 0$
points. But by Prop.  \ref{picard}, $X$ and $Y$ have the same Picard number,
so $r=0$, and $X$ and $Y$ are isomorphic.
\end{rmk}

\subsection{}
The following important result allows one to construct examples of FM
transforms. It was proved by Bondal and Orlov \cite{bo1}, and one of us \cite{br2}, using
ideas of Mukai.

\begin{thm}
\label{thea}
Suppose $X$ and $Y$ have dimension $n$.
Let $\P$ be an object of $\D(Y\times X)$, and let $\Phi$ denote the
exact functor
$$\Phi^{\P}_{Y\to X}\colon\D(Y)\lra\D(X)$$
defined
above. Then $\Phi$ is an equivalence of categories if and only if, for each point
$y\in Y$,
$$\Hom_{\D(X)}(\Phi(\OO_{y}),\Phi(\OO_{y}))=\C\mbox{ and }\Phi(\OO_y)\tensor\omega_X\cong\Phi\OO_y,$$
and for each pair of points $y_1,y_2\in Y$, and each integer $i$,
$$\Hom_{\D(X)}^i(\Phi(\OO_{y_1}),\Phi(\OO_{y_2}))=0\mbox{ unless }y_1=y_2\mbox{ and
}0\leq i\leq n.\qed$$
\end{thm}

Most examples of FM transforms for surfaces arise via the following
simple corollary. Recall that a sheaf on a smooth projective variety
$X$ is called \emph{special} if $E\tensor\omega_X= E$.

\begin{cor}
\label{fine}
Let $X$ be a smooth projective surface with a fixed polarisation, and
let $Y$ be a smooth, fine, complete, two-dimensional moduli space of
special, stable sheaves on $X$. Then there is a universal sheaf $\P$
on $Y\times X$, and the functor $\Phi^{\P}_{Y\to X}$ is a FM
transform.
\end{cor}

\begin{pf}
The assumption that $Y$ is fine means that there is a universal sheaf $\P$
on $Y\times X$, flat over $Y$. For each point $y\in Y$, $\P_y$ is a stable (hence
simple), special sheaf on $X$. Furthermore, since $Y$ is smooth of
dimension 2, the tangent space to $Y$ at $y$, which is given by
$$\Ext^1_X(\P_y,\P_y)$$
has dimension 2. It follows that for any pair of points $y_1, y_2\in
Y$,
$$\eu(\P_{y_1},\P_{y_2})=1-2+1=0.$$
If $y_1$ and $y_2$ are distinct then there are no non-zero maps
between the sheaves $\P_{y_1}$ and $\P_{y_2}$, so by Serre duality,
$$\Ext^2_X(\P_{y_1},\P_{y_2})=\Hom_X(\P_{y_2},\P_{y_1})^{\dual}=0.$$
The result then follows from Theorem \ref{thea}.
\qed
\end{pf}

\subsection{}
Assume that $X$ and $Y$ are surfaces. Our basic tool for classifying
FM transforms is

\begin{lemma}
\label{key}
Let $\Phi\colon\D(Y)\to \D(X)$ be a FM
transform and take a point $y\in Y$. Then there is an inequality
$$\sum_i\dim \Ext^1_X(\Phi^i(\OO_y),\Phi^i(\OO_y))\leq 2,$$
and moreover, each of the sheaves $\Phi^i(\OO_y)$ is special.
\end{lemma}

\begin{pf}
The second statement is immediate from Theorem \ref{thea}. For the
first part consider the spectral sequence \cite[Prop. 4.2]{bo1},
$$E^{p,q}_2=\bigoplus_i\Ext^p_X(\Phi^i(\OO_y),\Phi^{i+q}(\OO_y))\implies\Hom_{\D(X)}^{p+q}(\Phi(\OO_y),\Phi(\OO_y)).$$
The $E^{1,0}_2$ term survives to infinity, and by Lemma \ref{percy}
$$\Hom_{\D(X)}^1(\Phi(\OO_y),\Phi(\OO_y))=\Hom_{\D(Y)}^1(\OO_y,\OO_y)=\C^2,$$
so the result follows.
\qed
\end{pf}

\begin{cor}
\label{abelsur}
Suppose $X$ is an abelian surface. Then every FM transform $\Phi\colon\D(Y)\to\D(X)$ is a
sheaf transform.
\end{cor}

\begin{pf}
For any non-zero sheaf $E$ on an abelian surface $X$, the dimension of the space
$\Ext^1_X(E,E)$ is at least $2$.
\qed
\end{pf}

\subsection{}
The support of an object $E$ of $\D(X)$ is defined to be the
union of the supports of the cohomology sheaves of $E$. It is a closed
subset of $X$. A point $x\in X$ lies in the support of an object $E$
of $\D(X)$ if and only if there is an integer $i$ such that
$$\Hom_{\D(X)}^i(E,\OO_x)\neq 0.$$
This statement follows from a simple spectral sequence argument
\cite[Ex. 2.2]{br2}.

Suppose $\Phi\colon\D(Y)\to\D(X)$ is a FM transform, and let $\Psi$ be a
quasi-inverse $\D(X)\to\D(Y)$. Let $n$ be the common dimension of $X$ and
$Y$. For any pair of points $(y,x)\in Y\times X$,
\begin{eqnarray*}
\Hom_{\D(Y)}^i(\Psi(\OO_x),\OO_y)&=&
\Hom_{\D(X)}^i(\OO_x,\Phi(\OO_y)) \\
&=&\Hom_{\D(X)}^{n-i}(\Phi(\OO_y),\OO_x)^{\dual},
\end{eqnarray*}
so $x$ lies
in the support of $\Phi(\OO_y)$ precisely when $y$ lies in the
support of $\Psi(\OO_x)$.

A simple consequence of this is that the
supports of the objects $\Phi(\OO_y)$, as $y$ varies in $Y$, cover
$X$. For otherwise there would be a point $x\in X$ such that $\Psi(\OO_x)$
had empty support, and hence was zero, contradicting the assumption that $\Psi$ is an
equivalence. An extension of this argument gives

\begin{lemma}
\label{later}
Let $X$ and $Y$ be surfaces, and $\Phi\colon\D(Y)\to \D(X)$ a FM
transform. Suppose $X$ has non-zero Kodaira dimension, and take a
finite set of points $S\subset X$. Then for a general point $y\in Y$, the support of $\Phi(\OO_y)$ is disjoint
from S.
\end{lemma}

\begin{pf}Assume the opposite. Then
every point of $Y$ lies in the union over $x\in S$ of the supports of the
objects $\Psi(\OO_x)$, so for some $x\in S$, the support of $\Psi(\OO_x)$ is
the whole of $Y$. Since each cohomology sheaf of $\Psi(\OO_x)$ is
special, it follows that $\omega_Y$ has finite order, contradicting
Lemma \ref{extra}.
\qed
\end{pf}

%************************************************************************

\section{Ruled surfaces and surfaces of general type}

We start our classification of Fourier-Mukai transforms by looking at surfaces
with Kodaira dimension $-\infty$ and $2$.

\begin{prop}
\label{plop}
Let $X$ be a minimal surface of general type. Then the only FM partner
of $X$ is $X$ itself.
\end{prop}

\begin{pf}
It is a standard fact \cite[VII.2.3, VII.2.5]{bpvv} that
$X$ has only finitely many irreducible curves $D$ satisfying
$D\cdot \K_X= 0$. Thus, by Lemma \ref{later}, given a FM transform
$\Phi\colon\D(Y)\to\D(X)$,
we may choose $y\in Y$ so that the support of $\Phi(\OO_y)$ does not
contain any of these curves.

Let $E$ be a non-zero cohomology sheaf of
$\Phi(\OO_y)$. Since $E$ is special, $E$ is a torsion sheaf, and $\cl_1(E)\cdot \K_X=0$, so
$\cl_1(E)=0$, and $E$ is supported in dimension zero. Then Riemann-Roch gives $\eu(E,E)=0$, and this implies that $\Ext^1_X(E,E)$ has dimension at
least $2$. This applies to any cohomology sheaf of $\Phi(\OO_y)$ so Lemma \ref{key}
implies that some shift of $\Phi(\OO_y)$ is a sheaf $E$. Then, by Lemma
\ref{percy}, $E$ is simple, hence isomorphic to
$\OO_x$ for some $x\in X$, and Remark \ref{ratty} shows that $Y$ is isomorphic to $X$.
\qed
\end{pf}

\begin{prop}
\label{doob}
Let $X$ be a minimal surface of Kodaira dimension $-\infty$ with no
elliptic fibration. Then the only FM partner of $X$ is $X$ itself.
\end{prop}

\begin{pf}
Let $\Phi\colon\D(Y)\to\D(X)$ be a FM transform. If $X=\PP^2$ then the
argument of Prop. \ref{plop} applies, so we may
take $X$ to be a ruled surface over a curve of genus $g$.

We freely use notation and results from \cite[\S V.2]{ha2}.
Recall in particular that the N{\'e}ron-Severi group of $X$ is generated
by two elements $C$ and $f$, satisfying
$$C^2=-e, \quad f^2=0, \quad C\cdot f=1,$$
where $e$ is some integer invariant of $X$. In terms
of this basis
$$\K_X=-2C+(2g-2-e)f.$$

We shall assume for the moment that $X$ is not the unique rational ruled
surface with invariant $e=2$.

Suppose $D$ is an irreducible curve
on $X$ with $D\cdot\K_X=0$. We claim that $D^2\geq 0$.
Assume for contradiction that $D^2<0$.
The two-dimensional vector space
$\NS(X)\tensor_{\Z}\Q$ is then spanned by $D$ and $\K_X$, so the Hodge
index theorem implies that $\K_X^2>0$, and $X$ must be ruled over
$\PP^1$. In particular $e\geq 0$.
Write
$D=aC+bf$.
Since $D$ is irreducible
$$D\cdot C=b-ae\geq 0,\qquad D\cdot f=a\geq 0,$$
which is impossible since $D^2=a(2b-ea)<0$.

If $\Phi$ is not a sheaf transform then by Lemma \ref{key} we can
find a point $y\in Y$ and a cohomology sheaf $E$ of $\Phi(\OO_y)$,
supported in dimension 1,
such that the group $\Ext^1_X(E,E)$ has dimension at most 1. Since $E$
is special, any irreducible curve $D$ contained in the support of
$E$ satisfies $D\cdot\K_X=0$, and therefore $D^2\geq 0$.
But this is a contradiction since the group $\Ext^2_X(E,E)$ is non-zero,
so by Riemann-Roch, $\cl_1(E)^2<0$.

Thus the kernel of $\Phi$ is a sheaf $\P$ on $Y\times X$, flat over
$Y$. Given a point $y\in Y$ let $D$ be
the support of the sheaf
$\P_y=\Phi(\OO_y)$. If $D$ is zero-dimensional, it follows as in Prop.
\ref{plop} that $X$ and $Y$ are isomorphic. Thus we may assume that
$D$ is  a curve. Then
$D\cdot\K_X=0$ and $D^2=0$, and since
$\P_y$ is simple, $D$ is connected. Furthermore $D$ is
irreducible since any irreducible component $D_0\subset D$ satisfies
$D_0\cdot\K_X=0$ and hence $D_0^2\geq 0$.

Fix a smooth hyperplane section $H$ on $Y$, with $H\cdot\K_Y\neq
0$, and let $\Psi\colon\D(X)\to\D(Y)$
be a quasi-inverse of $\Phi$.
For any point $x\in X$, the
support of $\Psi(\OO_x)$ meets $H$ at a
finite set of points because $\Psi(\OO_x)$ is special. We show that
for some integer $d$ this defines a map $X\to \Sym^d(H)$
which is an elliptic fibration on $X$.

Recall the definition of  the
derived dual
$$\OO_H^{\dual}=\R\lHom_{\OO_Y}(\OO_H,\OO_Y)=\OO_H(H)[-1].$$
For any line bundle
$L$ on $Y$, one has isomorphisms
$$\Hom^i_{\D(Y)}(L,\Psi(\OO_x)\Ltensor\OO_H)=\Hom^i_{\D(Y)}(L\tensor\OO_H^{\dual},\Psi(\OO_x))$$
$$=\Hom^{i+1}_{\D(X)}(\Phi(L\tensor\OO_H(H)),\OO_x).$$
By the theorem on cohomology and base-change, we can choose $L$
sufficiently ample so that the object
$$\Phi(L\tensor\OO_H(H))=\R\pi_{X,*}(\P\tensor \pi_Y^*(L\tensor\OO_H(H)))$$
is concentrated in degree 0 and is locally free. Then
the above groups vanish unless $i=-1$, so for each $x\in X$, the
object
$$\Psi(\OO_x)[-1]|_H=\Psi(\OO_x)[-1]\Ltensor\OO_H$$
is concentrated in degree 0. If the kernel of the
transform $\Psi[-1]$ is the object $\QQ$ of $\D(Y\times X)$, this implies
\cite[Lemma 4.3]{br2}
that
$\QQ|_{H\times X}$
is a sheaf on $H\times X$, and is
flat over $X$. Thus
$\QQ$ defines a family of torsion sheaves on
$H$, parameterised by $X$, so induces a morphism
$$f\colon X\to\Sym^d(H).$$

By the result of Section 2.7, any fibre of $f$ is
the intersection over a finite set of points $y\in H$ of the supports
of the sheaves  $\P_y=\Phi(\OO_y)$. Each sheaf
$\P_y$ is supported on
an irreducible curve $D_y$ satisfying $D_y\cdot\K_X=0$, and by
Riemann-Roch, given two points $y_1,y_2\in Y$ one has $D_{y_1}\cdot
D_{y_2}=0$.
It follows that any non-singular fibre of $f$ is an elliptic
curve.
Applying Stein factorisation gives an
elliptic fibration $X\to S$ onto a smooth curve $S$,
and hence a contradiction.

The remaining possibility is that $X$ is the unique rational ruled
surface with invariant $e=2$. Then, \cite[Cor. 2.18]{ha2},
$C\subset X$ is the only irreducible curve
satisfying $C\cdot\K_X=0$, so the argument of Prop. \ref{plop}
shows that $Y$ is birational to $X$.
By Lemma
\ref{picard}, $X$ and $Y$ have the same Picard number, so $Y$ is also
a rational ruled surface, and hence has no elliptic fibration.
Applying what we have already proved to
$Y$, we conclude that $Y$ also has
invariant $e=2$, so $X$ and $Y$ are isomorphic.
\qed
\end{pf}

%************************************************************************

\section{Elliptic surfaces}

Fourier-Mukai transforms for elliptic surfaces were introduced in \cite{br1}. We start by reviewing
the
construction given there. Throughout we fix a surface $X$ and a
relatively minimal elliptic fibration $\pi\colon X\to C$.

Given an
object $E$ of $\D(X)$, one defines the \emph{fibre degree} of $E$
$$\dg(E)=\cl_1(E)\cdot f,$$
where $f$ denotes the algebraic equivalence class of a fibre of $\pi$.
Let $\m_{X/C}$
denote the highest
common factor of the fibre degrees of objects of
$\D(X)$. Equivalently, $\m_{X/C}$ is the smallest number $d$ such that there is a holomorphic
$d$-section of $\pi$.

Let $a>0$ and $b$ be integers, with $b$ coprime to
$a\m_{X/C}$. Then, as was shown in \cite{br1}, there is a smooth,
two-dimensional component
$$Y=\J_{X/C}(a,b)$$ of the moduli
space of pure dimension one stable sheaves on $X$, the general point of
which represents a rank $a$, degree $b$ stable vector bundle
supported on a smooth fibre of $\pi$.

There is a natural morphism
$Y\to C$, taking a point representing a sheaf supported
on the fibre $\pi^{-1}(p)$ of $X$ to the point $p$, and this morphism is a
relatively minimal
elliptic fibration. Moreover, there is a universal sheaf $\P$ on
$Y\times X$, supported on $Y\times_C X$, and the corresponding functor
$\Phi^{\P}_{Y\to X}$ is a FM transform.
In \cite{br1} these ideas
are used to prove the following result.

\begin{thm}
\label{super}
Take an element
$$\mat{c}{a}{d}{b}\in\mbox{SL}_2(\Z),$$
such that $\m_{X/C}$ divides $d$ and $a>0$. Let $Y$ be the elliptic
surface $\J_{X/C}(a,b)$ over $C$. Then there exist universal sheaves $\P$ on
$Y\times X$, flat over both factors, such that the resulting functor $\Phi=\Phi^{\P}_{Y\to X}$ is an
equivalence of categories satisfying
\begin{equation}
\label{num}
\col{\rk(\Phi E)}{\dg(\Phi E)}=\mat{c}{a}{d}{b}\col{\rk(E)}{\dg(E)},
\end{equation}
for all objects $E$ of $\D(Y)$.
\qed
\end{thm}

When $a=1$ the elliptic surface $\J_{X/C}(a,b)$ is the relative Picard scheme
of R. Friedman \cite{fr}, which we denote more concisely by $\J_{X/C}(b)$.

\begin{lemma}
\label{dave}
For any pair of integers
$a$, $b$, with $b$ coprime to $a\m_{X/C}$, there is an isomorphism
\begin{equation}
\label{woops}
\J_{X/C}(a,b)\isom \J_{X/C}(b).
\end{equation}
\end{lemma}

\begin{pf}
The essential point is that, as in \cite[Theorem 6]{at},  the determinant map gives an isomorphism
$$\det\colon\M_{X_p}(a,b)\lra \M_{X_p}(1,b)=\Jac^b(X_p)$$
on each smooth fibre $X_p$ of $\pi\colon X\to C$.
Glueing these maps gives the isomorphism (\ref{woops}).

In more detail, let $U\subset C$ be an open subset
of $C$ over which the morphism $\pi$ is smooth, and let $\P_U$ denote
the restriction of $\P$ to the open subset
$$Y_U\times_C X_U\subset Y\times_C X.$$
For each point $y\in Y_U$, the restriction of $\P_U$ to the fibre
$\{y\}\times_C X$ is a rank $a$, degree $b$ vector bundle on the
smooth elliptic curve $X_{\pi(y)}$. Therefore $\P_U$ is
locally free, and we can consider the determinant line bundle
$\wedge^a \,\P_U$.
This parameterises degree $b$ line bundles on the fibres of $\pi$, and
hence defines an isomorphism
$$\J_{X/C}(a,b)\times_C U\lra\J_{X/C}(b)\times_C U.$$
Since both spaces in (\ref{woops}) are relatively minimal over $C$,
\cite[Prop. III.8.4]{bpvv} implies that they are isomorphic.
\qed
\end{pf}

\begin{lemma}
For any integer $b$ coprime to $\m_{X/C}$, the elliptic surface
$Y=\J_{X/C}(b)$ has the same Kodaira dimension as $X$.
\end{lemma}

\begin{pf}
The Euler
numbers of $X$ and $Y$ are equal by Lemma \ref{picard}. By
\cite[Prop. I.3.23]{frmo} we must show that the base orbifolds of $X$
and $Y$ are diffeomorphic, i.e. that for each point $p\in C$, the
multiplicities of the fibres of $X$ and $Y$ over the point $p$ are
equal.

Let $\Phi\colon\D(Y)\to\D(X)$ be a FM transform as in Theorem
\ref{super}, and fix a point $p\in C$. Let $F_Y$ be the fibre of $Y$ over $p$. Thus $\OO_{F_Y}$ has Chern character
$(0,f,0)$. It is easy to check using formula (\ref{num}) that the object $\Phi(\OO_{F_Y})$ of
$\D(X)$ has Chern
character $-(0,cf,d)$.

Suppose that $F_Y$ is a multiple fibre, so that $F_Y=mD_Y$ for some positive integer $m$, and let
$E=\Phi(\OO_{D_Y})$. Then $E$ is an object of $\D(X)$ supported on the
fibre $F_X$ of $X$ over $p$, and if the first Chern class of $E$ is $-D_X$,
we must have $mD_X=cf$. But $c$ is coprime to $d$, and $m$ divides
$\m_{X/C}$ which divides $d$, so it follows that $F_X$  has multiplicity at least $m$. By symmetry the multiplicities of
the fibres of $X$ and $Y$ over $p$ are equal.
\qed
\end{pf}

We can now prove the following classification result.

\begin{prop}
\label{elika}
Let $X$ be a minimal surface of non-zero Kodaira dimension, with an
elliptic fibration $\pi\colon X\to
C$. If $Y$ is an FM partner of $X$,
then $Y$ is isomorphic to the relative Picard scheme $\J_{X/C}(b)$, for some
integer $b$ coprime to $\m_{X/C}$.
\end{prop}

\begin{pf}
Let $f$ be the cohomology class of a fibre of $\pi$. The cohomology class $\K_X$
is a non-zero rational multiple of $f$, \cite[Cor. V.12.3]{bpvv}, so the support of any special
sheaf on $X$ is contained in a finite number of fibres of $\pi$.

Take $x\in X$ lying on a
smooth fibre of $\pi$, and take a point $y\in Y$ such that the support
of the object
$E=\Phi(\OO_y)$ contains $x$. Since $\Hom_{\D(X)}(E,E)=\C$, the
support of $E$ is connected, hence equal to the (smooth) fibre of $\pi$
containing $x$. Now the Chern class
of $E$ must be $(0,af,-b)$ for some integers $a$ and $b$, and since
$$\eu(E,\Phi(\OO_Y))=\eu(\OO_y,\OO_Y)=1,$$
Riemann-Roch implies that $a\m_{X/C}$ is coprime to $b$. Since $E$ is
supported on an elliptic curve, all of its cohomology sheaves are non-rigid, so Lemma \ref{key} implies that $E$ has only one
non-zero cohomology sheaf. Thus a shift of $E$ is a
simple sheaf on an elliptic curve, hence stable.

Let $Y^+$ be
the elliptic surface $\J_{X/C}(b)$, with its relatively minimal elliptic
fibration $\pi^+\colon Y^+\to C$. There is a
transform
$$\Psi\colon\D(Y^+)\lra \D(X)$$
 which takes the structure sheaf of
some point of $Y^+$ to $E$. Applying Prop. \ref{birat} to the transform
$\Psi^{-1}\circ\Phi$ shows
that there is a birational equivalence $f\colon Y\dashrightarrow Y^+$,
such that
$$(\Psi^{-1}\circ\Phi)(\OO_y)=\OO_{f(y)}$$
for all points $y$ in some open subset of $Y$.

If $X$ has Kodaira
dimension 1, then so do $Y$ and $Y^+$, so Remark \ref{ratty} shows
that $f$ extends to an isomorphism, completing the proof. The only
other possibility is that $X$ is a minimal ruled surface over an elliptic
base. In that case $Y$ and $Y^+$ also have Kodaira dimension
$-\infty$, and also have Picard number $2$, so are minimal ruled
surfaces over an elliptic base. By Prop. \ref{doob} we may assume that $Y$ has an
elliptic fibration $\pi\colon Y\to C$.

Consider the full subcategory $\D_{sp}(Y)$ of $\D(Y)$
consisting of objects invariant under twisting by $\omega_Y$. The
support of any object of $\D_{sp}(Y)$ is contained in the union of a
finite number of fibres of $\pi$. By uniqueness of Serre functors (see
Lemma \ref{extra}), the FM transform $\Psi^{-1}\circ\Phi$ takes objects of
$\D_{sp}(Y)$ to objects of $\D_{sp}(Y^+)$. This says that the
birational map $f$ takes fibres of $\pi$ to fibres of $\pi^+$, so
applying \cite[Prop. III.8.4]{bpvv} shows that $f$ extends to an isomorphism.
\qed
\end{pf}

\begin{rmk}
\label{misc}
If $\sigma$ is a divisor on $X$ such that $\sigma\cdot
f=\m_{X/C}$, then twisting by $\OO_X(\sigma)$ gives an isomorphism
$$\J_{X/C}(b)\isom\J_{X/C}(b+\m_{X/C}).$$
Thus an elliptic surface of non-zero Kodaira dimension has only
finitely many FM partners.

The argument of Lemma \ref{dave} shows that the operation of taking duals
gives a
birational equivalence, hence an isomorphism
$$\J_{X/C}(b)\isom\J_{X/C}(-b).$$

Finally note that there is an isomorphism
$$\J_{X/C}(1)\isom X.$$
To see this note that the ideal sheaf of the diagonal in $X\times_C X$ is
flat over both factors, and hence generates a family of rank 1, degree
$-1$ stable sheaves supported on the fibres of $\pi$.
\end{rmk}

% ***************************************************************************

\section{K3 and abelian surfaces}

Let $X$
be an abelian or K3 surface. Following Mukai \cite{muk5}, one introduces the
\emph{extended Hodge lattice} of $X$ by using the formula
$$\lb(r_1,D_1,s_1),(r_2,D_2,s_2)\rb=D_1\cdot D_2-r_1 s_2-r_2 s_1$$
to define a bilinear form on the
cohomology ring
$$\H^{\ts}(X,\Z)=\H^0(X,\Z)\oplus \H^2(X,\Z)\oplus \H^4(X,\Z),$$
and taking the following Hodge decomposition of $\H^{\ts}(X,\C)$
$$\H^{\ts\,(0,2)}(X,\C)=\H^{0,2}(X,\C),\qquad
\H^{\ts\,(2,0)}(X,\C)=\H^{2,0}(X,\C),$$
$$\H^{\ts\,(1,1)}(X,\C)=\H^0(X,\C)\oplus \H^{1,1}(X,\C)\oplus \H^4(X,\C).$$

Inside $\H^2(X,\Z)$ one has two sublattices, the \emph{N{\'e}ron-Severi
lattice} which is
$$\NS(X)=\H^2(X,\Z)\cap \H^{1,1}(X,\C),$$
and its orthogonal complement $\T(X)$, the \emph{transcendental
lattice} of $X$. The transcendental lattice inherits a Hodge structure
from $\H^2(X,\Z)$.

A Hodge isometry $\phi$ between the transcendental lattices (or
extended Hodge lattices) of two K3 (or abelian) surfaces $X$ and $Y$, is
an isometry between the relevant lattices which preserves the Hodge
decomposition. This last condition is equivalent to the
statement that $\phi\tensor\C$ takes the cohomology class of the unique (up to scalar
multiples) non-vanishing holomorphic 2-form on $Y$ to the
corresponding class on $X$.

To each sheaf (or complex of sheaves) $E$ on $X$ one associates a
\emph{Mukai vector}
$$v(E)=(\rk(E),\cl_1(E),\ha \cl_1(E)^2-\cl_2(E)-\epsilon \rk(E))\in \H^{\ts}(X,\Z),$$
where $(\rk(E),\cl_1(E),\cl_2(E))$ are the Chern classes of $E$, and
$\epsilon$ is 0 or 1 depending on whether $X$ is abelian or K3
respectively. Having done this, the Riemann-Roch formula becomes
$$\eu(E,F)=-\lb v(E),v(F) \rb,$$
for any pair of sheaves (or complexes) $E$ and $F$ on $X$.

By Lemma \ref{extra} any FM partner of
$X$ is also an abelian or K3 surface, and Lemma \ref{picard} shows that an abelian surface could never be a partner of a
K3 surface.

The following theorem is due to Mukai \cite{muk5} and Orlov \cite{orl1}. We sketch the proof for the reader's convenience, and to fix
ideas for the next section where similar techniques are used.

\begin{thm}
Let $X$ and $Y$ be a pair of K3 (respectively abelian) surfaces. The following statements are equivalent,
\begin{itemize}
\item[(a)] there is a FM transform $\Phi\colon\D(Y)\to\D(X)$,

\item[(b)] there is a Hodge isometry $\phi^t\colon\T(Y)\to\T(X)$,

\item[(c)] there is a Hodge isometry
$\phi\colon\H^{\ts}(Y,\Z)\to\H^{\ts}(X,\Z)$,

\item[(d)] $Y$ is isomorphic to a fine, two-dimensional moduli space of
stable sheaves on $X$.
\end{itemize}
\end{thm}

\par
\noindent {\it Sketch proof.}

(a)$\implies$ (b). Any FM transform $\Phi\colon\D(Y)\to\D(X)$ induces an
isomorphism of vector spaces
$$\phi\colon \H^{\ts}(Y,\C)\lra \H^{\ts}(X,\C),$$
as in Section \ref{cc}. Since the kernel of $\Phi$ is
algebraic, this isomorphism preserves the Hodge decomposition. Mukai checks
\cite[Lemma 4.7, Theorem 4.9]{muk5},
that $\phi$ preserves the inner product and the integral
lattices. It follows that $\phi$ takes $\T(Y)$ into $\T(X)$.

\smallskip

(b)$\implies$ (c). This is a consequence of a result of Nikulin
\cite[Prop. 6.1]{muk5}. The orthogonal complement of $\T(Y)$ in $\H^{\ts}(Y,\Z)$ contains the
hyperbolic lattice
$$H=\H^0(Y,\Z)\oplus \H^4(Y,\Z),$$
 so any isometry of
transcendental lattices extends to an isometry of extended Hodge lattices.

\smallskip

(c)$\implies$ (d). Let
$$\phi\colon \H^{\ts}(Y,\Z)\lra \H^{\ts}(X,\Z)$$
be a Hodge isometry, and put
$v=\phi(0,0,1)$. Composing with standard automorphisms of
$\H^{\ts}(X,\Z)$, obtained either by swapping $H^0$ and $\H^4$, or by twisting
by line bundles, we may assume that
$v=(r,\ell,s),$
 with $r$ positive, $\ell$
ample, and $s$ coprime to $r$.

Since $v$
is algebraic, we can consider $Y^+$, the moduli space of stable sheaves
on $X$ with Mukai vector $v$, with respect to the polarization $\ell$.
The fact that $v$ is primitive implies
that this moduli space is fine \cite[Theorem A.6]{muk5}, and the fact that $v^2=0$ implies
that $Y^+$ is two-dimensional. General results of Mukai show that $Y^+$
is smooth \cite{muk3} and non-empty \cite[Theorem 5.4]{muk5},
\cite[Prop. 6.16, Cor. 6.23]{muk1}.

By Cor. \ref{fine} there is a FM transform
$$\Psi\colon\D(Y^+)\lra\D(X),$$
such that for any point $y\in Y^+$,
$\Psi(\OO_y)$
is the corresponding stable sheaf on $X$. The argument given for
(a)$\implies$(b) shows that $\Psi$ gives rise to a Hodge
isometry
$$\psi\colon \H^{\ts}(Y^+,\Z)\lra \H^{\ts}(X,\Z)$$
taking $(0,0,1)$ to $v$. The composite $\psi^{-1}\circ\phi$
is a Hodge isometry
$$\H^{\ts}(Y,\Z)\lra \H^{\ts}(Y^+,\Z)$$
fixing
$(0,0,1)$, which therefore restricts to give a Hodge isometry
$$\H^2(Y,\Z)\lra
\H^2(Y^+,\Z).$$

In the K3 surface case, the Torelli theorem shows that $Y$ and
$Y^+$ are isomorphic, and we are done. Otherwise $X$ is an abelian
surface, and \cite[Theorem 1]{shi} shows
that  $Y$ is isomorphic to either $Y^+$ or its dual variety. In either case $Y$ is a
FM partner of $X$ since dual abelian varieties have equivalent derived
categories by the
results of \cite{muk2}. It follows from  Cor. \ref{abelsur} that there is a
universal family of sheaves
$\{\P_y:y\in  Y\}$ on $X$, which we may assume are locally free, and which are simple by Lemma \ref{percy}. Then \cite[Prop. 6.16]{muk1} shows that
each bundle $\P_y$ is actually stable (with respect
to any polarization of $X$),
so $Y$ is indeed a moduli space of stable
sheaves on $X$.

\smallskip

(d)$\implies$ (a). Immediate from Cor. \ref{fine}.
\qed

\begin{rmk}
Given a FM transform $\Phi\colon\D(Y)\to \D(X)$ between K3 surfaces, the
theorem implies that $Y$ is a moduli space of stable sheaves on
$X$. This does \emph{not} mean that $\Phi$ is given by the formula
$\Phi^{\P}_{Y\to X}$, with $\P$ a universal sheaf on $Y\times X$. As we
mentioned in the introduction, finding the set of FM transforms between two
K3 surfaces satisfying the conditions of the theorem is a difficult
unsolved problem.
\end{rmk}

\begin{prop}
Let $X$ be a K3 or abelian surface. Then $X$ has only a finite number of FM
partners.
\end{prop}

\begin{pf}
Suppose for contradiction that
$X$ has infinitely many FM partners $Y$, and choose two such surfaces
$Y_1$ and $Y_2$,
not isomorphic,
together with a Hodge isometry
$$\phi\colon \H^{\ts}(Y_1,\Z)\lra \H^{\ts}(Y_2,\Z).$$
Since $\phi$ preserves the Hodge decomposition, it induces lattice isomorphisms
$$\T(Y_1)\lra\T(Y_2), \qquad \NS(Y_1)\oplus H\lra \NS(Y_2)\oplus H,$$
where $H=\H^0(Y,\Z)\oplus \H^4(Y,\Z)$ is the hyperbolic lattice. The second isomorphism shows that
the lattices $\NS(Y_1)$ and $\NS(Y_2)$ have the same genus
\cite[Theorem 1.3.1, Cor. 1.9.4]{nik}. There are only
finitely many different lattices of each genus \cite[Ch. 10, \S
3.3]{cafr}, so we may choose $Y_1$ so that there are infinitely many
possible choices for $Y_2$ such that $\NS(Y_1)$ and $\NS(Y_2)$
are isometric. For any such choice we can find a Hodge isometry
$$f\colon\NS(Y_1)\oplus \T(Y_1)\lra \NS(Y_2)\oplus \T(Y_2).$$

Fix an abstract lattice $W$
isomorphic to
$\NS(Y_1)\oplus\T(Y_1)$.
Lattices containing $W$ as a sublattice of finite
index are all contained in the dual lattice $W^*=\Hom_{\Z}(W,\Z)$, and thus correspond to subgroups of
the finite abelian group $W^*/W$, as in \cite[\S4]{nik}. Obviously there are only a finite
number of these, so, changing $Y_1$ again, we may assume that there
are infinitely many possible choices for $Y_2$ such that the lattice
extensions
$$\NS(Y_1)\oplus\T(Y_1)\embed \H^2(Y_1,\Z),\qquad \NS(Y_2)\oplus\T(Y_2)\embed \H^2(Y_2,\Z)$$
are isomorphic. But for any such choice, the isometry $f$ extends to a Hodge isometry
$$\H^2(Y_1,\Z)\lra \H^2(Y_2,\Z).$$
If $X$ is a K3 surface, the Torelli theorem implies that $Y_1$ and $Y_2$ are isomorphic. In the case when $X$ is an abelian surface we can
apply \cite[Theorem 1]{shi} to conclude that $Y_1$ is isomorphic to $Y_2$
or its dual. In both cases we obtain a contradiction, since we claimed
there were an infinite number of possible choices for $Y_2$.
\qed
\end{pf}

\begin{rmk}
Nikulin's results imply that if a K3 surface has
Picard number at least 12 then it has no FM partners other than
itself \cite[Prop. 6.2]{muk5}. This result is not true in general;
Mukai observes \cite[p. 394]{muk5} that there are K3 surfaces with isometric
transcendental lattices (hence equivalent derived categories) but
distinct N{\'e}ron-Severi lattices.
\end{rmk}

% ***************************************************************************

\begin{section}{Enriques and bielliptic surfaces}

We conclude our classification of FM transforms by considering surfaces
with non-trivial canonical bundle of finite order, i.e. Enriques and
bielliptic surfaces. Collectively we shall call such surfaces
\emph{quotient surfaces}.

\smallskip

All bielliptic surfaces have exactly two elliptic fibrations, and the
general Enriques surface is also an elliptic surface in two different
ways. Thus it is possible to generate many examples of non-trivial FM
transforms for quotient surfaces by considering compositions of the transforms
arising from Theorem
\ref{super}. Further examples were described in \cite{br3}. Nonetheless, in this
section we shall prove that if $X$ is a quotient surface then any FM partner of $X$ is
isomorphic to $X$ itself.

\smallskip

Let $X$ be a quotient surface, and let $n$ be the order of
$\omega_X$. It is easily seen that there is a surface
$\Xt$, with trivial canonical bundle, such that $X$ is the quotient of
$\Xt$ by a free action of the finite cyclic group $G$ of order $n$. We
refer to the quotient map
$p_X\colon\Xt\to X$ as the \emph{canonical cover} of $X$.
Let $$\Phi\colon\D(Y)\to\D(X)$$ be a FM
transform. By Lemma \ref{extra}, $\omega_Y$ also has order $n$,
so has a canonical cover $p_Y\colon\Yt\to Y$. In \cite{br3} we proved that there
is a \emph{lift} of $\Phi$ to a FM transform
$$\Phit\colon\D(\Yt)\lra\D(\Xt)$$
making the following two squares of functors commute
$$
\begin{array}{ccc}
\D(\Yt) &\lRa{\Phit} & \D(\Xt) \\
\scriptstyle{p_Y^*}\bua\bda\scriptstyle{p_{Y,*}}
&&\scriptstyle{p_X^*}\bua\bda\scriptstyle{p_{X,*}} \\
\D(Y) &\lRa{\Phi} &\D(X).
\end{array}
$$
Moreover any such lift is equivariant, in that there is an automorphism
$\mu\colon G\to G$ such that for each $g\in G$ there is an isomorphism
of functors
$$g^*\circ\Phit\isom \Phit\circ\mu(g)^*.$$

\begin{prop}
Let $X$ be an Enriques surface. Then the only FM partner of
$X$ is $X$ itself.
\end{prop}

\begin{pf}
Take notation as above. Thus $\Xt$ and $\Yt$ are K3 surfaces and $G=\Z/(2)$. It
follows that $Y$ is also an Enriques surface.
Let $g$ be the generator of $G$. Define sublattices
$$\H^{\ts}_{\pm}(\Xt,\Z)=\{\theta\in\H^{\ts}(\Xt,\Z):g^*(\theta)=\pm\theta\}.$$
Note that $\H^{\ts}_{-}(\Xt,\Z)\subset\H^2(\Xt,\Z)$. Furthermore
$$\H^{0,2}(\Xt,\C)\subset \H^{\ts}_{-}(\Xt,\Z)\tensor\C.$$

The transform $\Phit$ induces
a $G$-equivariant Hodge isometry between $\H^{\ts}(\Yt,\Z)$ and
$\H^{\ts}(\Xt,\Z)$, and hence gives
a $G$-equivariant isometry
$$f\colon \H^{\ts}_{-}(\Yt,\Z)\lra \H^{\ts}_{-}(\Xt,\Z),$$
taking the subspace $\H^{0,2}(\Yt,\C)$
onto $\H^{0,2}(\Xt,\C)$. We claim that $f$ extends to an isometry
$$f\colon\H^2(\Yt,\Z)\lra \H^2(\Xt,\Z).$$
Assuming this for the moment, note that $f$ is then a
$G$-equivariant Hodge isometry, so by the Torelli theorem for Enriques
surfaces,
\cite[VIII.21.2]{bpvv}, $X$ and $Y$ are isomorphic.

To prove the claim we use more results of Nikulin.
The orthogonal complement of $\H^{\ts}_{-}(\Xt,\Z)$ in
$\H^2(\Xt,\Z)$, which is equal to $\H^2_{+}(\Xt,\Z)$, is
even, 2-elementary (\cite[Defn. 3.6.1]{nik}) and indefinite.
The claim then follows from
Prop. 1.14.1 and Theorems 3.6.2 and 3.6.3 of \cite{nik}.
\qed
\end{pf}

\begin{prop}
\label{a}
Let $X$ be a bielliptic surface. Then the only FM partner of $X$ is
$X$ itself.
\end{prop}

\begin{pf}
Take notation as above. Then
$\Xt$ is
a quotient of a product of two elliptic curves by a
finite group $H$ of translations, so is an abelian surface
with two elliptic fibrations without
multiple fibres. The N{\'e}ron-Severi group of $\Xt$ is generated by the
algebraic equivalence classes $f_1,f_2$ of the corresponding fibres\footnote{This statement is false in general, although fortunately Proposition \ref{a} is nonetheless correct. See Section \ref{err}.} 
and $f_1\cdot f_2=m$, the order of $H$.

Note that the group $G$ acts on one of the fibres, say $f_2$, of $\Xt$ via
multiplication by a complex $n$th root of unity. It follows that $X$
has a multiple fibre of multiplicity $n$, and that there exists a
divisor $D$ on $X$ such that $p_X^* D=f_1$.

Consulting the table on
\cite[p. 148]{bpvv}, the possible values of $m$ are
1,2,3 and 4, and when $m>1$, the prime factors of $m$ are the
same as those of $n$. By Remark \ref{misc}, all the relative Picard
schemes of
$\Xt$ considered in Section 4  are isomorphic to $\Xt$. We
shall show that $\Phit$ is isomorphic to a composite of transforms
arising from the two elliptic fibrations via Theorem \ref{super}.

Cor. \ref{abelsur} shows that $\Phit$ is a sheaf transform. Thus
we can
suppose that for each point $\yt\in\Yt$, the object $\Ft=\Phit(\OO_{\yt})$ is
a sheaf on $\Xt$, of Chern character $(r,pf_1+qf_2,s)$ say.
The fact that $\Phit$ is a lift of $\Phi$ implies that
$$\eu(p_X^*(\Phi(\OO_Y)),\Ft)=\eu(\Phi(\OO_Y),p_{X,*}(\Ft))=\eu(\Phi(\OO_Y),\Phi(\OO_y))=1,$$
where $y=p_Y(\yt)$. Now $p_X^*(\Phi(\OO_Y))$ has second Chern class
divisible by $n$, and hence by $m$,
so Riemann-Roch implies that $s$ is coprime
to $m$.

Let $h$ be the greatest common divisor of $r$ and $p$, and take
integers $b$ and $d$ such that
$bp+dr=-h$.
Riemann-Roch together with Lemma \ref{percy} shows that $rs=mpq$. Since
$s$ is coprime to $m$, we see that $m$ divides
$r/h$. By Theorem \ref{super} we can find a transform $\Psit\colon\D(\Xt)\to\D(\Xt)$
such that
$$
\col{\rk(\Psit E)}{\dg(\Psit E)}=\mat{-p/h}{r/mh}{dm}{b}\col{\rk(E)}{\dg(E)},
$$
where $\dg(E)=\cl_1(E)\cdot f_2$ for any object $E$ of
$\D(\Xt)$.

The transform $\Psit$ arises by considering the moduli
space of stable sheaves $\Psit(\OO_{\yt})$ on $\Xt$ of Chern character $(0,af_2,b)$, where
$a=r/mh$. We must show that $\Psit$ is the lift of a transform
$\Psi\colon\D(X)\to\D(X)$, it follows from this that $\Psit$ is $G$-equivariant. By \cite[Lemma 5.1]{br3}, it is enough to
check that for some object $A$ of $\D(X)$,
$$\eu(p_X^* A,\Psit(\OO_{\yt}))=1.$$
By Riemann-Roch, this is the statement that there is a divisor $D$ on
$X$ such that $b$ is coprime to $p_X^* D\cdot f_2$. But we can assume
that $p_X^* D=f_1$, and this is enough, since $b$ is coprime to $m$.

Replacing $\Phi$ with the composite transform $\Psi\circ\Phi$ we can
now assume that $r=0$ and $p\neq 0$. By Riemann-Roch, $q=0$
also, so $\Ft$ has Chern character $(0,pf_1,s)$, where, as before, $s$ is coprime to
$p$ and $m$. There exists an equivariant transform
$$\Psit\colon\D(\Xt)\lra\D(\Xt)$$
such that
$\Psi(\OO_{\xt})$ has this same Chern character, so composing $\Phit$ with
the inverse of $\Psit$ we can assume that $\Ft=\Phit(\OO_{\yt})$ has Chern character
$(0,0,1)$. But $\Phit$ is a sheaf transform, so there is an isomorphism
$\tilde{\phi}\colon\Yt\to \Xt$ such that for all $\yt\in\Yt$,
$$\Phit(\OO_{\yt})=\OO_{\tilde{\phi}(\yt)}.$$
Since $\Phit$ is $G$-equivariant, $\tilde{\phi}$ descends to an
isomorphism $\phi\colon Y\to X$.
\qed
\end{pf}

\end{section}

\section{Erratum}
\label{err}

The proof of Proposition 6.2 in the published version of this paper is incorrect, since it assumes that the N{\'e}ron-Severi group of the canonical cover always has rank 2. We would like to thank  Rory Potter and Evgeny Shinder for pointing this out. Fortunately Proposition 6.2 is nonetheless correct. We now explain the proof, much of which is due to Potter\footnote{R.D. Potter, Derived autoequivalences of bielliptic surfaces, arxiv 1701.01015; R.D. Potter,  Derived categories of surfaces and group actions, PhD thesis, University of Sheffield, 2017, available at http://etheses.whiterose.ac.uk/19643/.}. Throughout, we denote the group of numerical equivalence classes of divisors on a smooth surface $S$ by $\Num(S)$. 

\smallskip

Consider a bielliptic surface $X$ with canonical bundle $\omega_X$ of order $n>1$. As in \cite[Section V.5]{bpvv}, the surface $X$ is a quotient $(A\times B)/J$, with $A$ and $B$ elliptic curves, and $J$ a finite subgroup of the translation group of $A$ which acts faithfully on $B$ by automorphisms. There is a decomposition $J=G\oplus H$, with $G,H\subset J$ both cyclic subgroups, such that $G$ has order $n$, and $H$ acts  on $B$ by translations. Writing $k$ for the order of $H$, the possible values of  $(n,k)$ are
\[(2,1), \quad (3,1), \quad (4,1), \quad (6,1),\quad (2,2), \quad (3,3), \quad (4,2).\]

The group $H$ acts on $A\times B$ entirely by translations, and the canonical cover $\tilde{X}$ is the quotient $(A\times B)/H$, and is an abelian surface.
The surface $\tilde{X}$ has two smooth elliptic fibrations \[\tilde{\pi}_A\colon \tilde{X}\to A/H, \qquad \tilde{\pi}_{B}\colon \tilde{X}\to B/H.\] We denote the classes in $\Num(\tilde{X})$ defined by the fibres of these maps by $\tilde{B}$ and $\tilde{A}$ respectively.
Similarly there are elliptic fibrations
\[{\pi}_A\colon {X}\to A/J, \qquad {\pi}_{B}\colon {X}\to B/J,\]
whose generic fibres define classes $B$ and $A$ in $\Num(X)$ respectively. %The first of these is a smooth fibration, but the second has multiple fibres.

Since $H^2(X,\OO_X)=0$,  the N{\'e}ron-Severi group $\NS(X)$  coincides with $H^2(X,\Z)$, and the numerical equivalence group $\Num(X)$ is the quotient of $\NS(X)$ by the torsion subgroup.
 Serrano proved\footnote{Theorem 1.4 in F. Serrano, Divisors of bielliptic surfaces and embeddings in $\mathbb P^4$, Math. Z., 203, 527 -- 533, (1990).}  that there are classes $A',B'\in \Num(X)$ with $nA'=A$ and $kB'=B$ such that \[\Num(X)=\Z A'\oplus \Z B'.\] 
Since $A\cdot B=kn$, it follows that $A'\cdot B'=1$. 

The canonical cover $p=p_X\colon \tilde{X}\to X$  induces a map
\[p_{*}\colon \Num(\tilde{X})\to \Num(X).\]
Consider the image  $\Delta=\operatorname{im}(p_*)\subset \Num(X)$ of this map. We claim that given $d\in \Z$ there are implications
\begin{equation}
\label{err_stuff}dA'\in \Delta\iff n\mid d, \qquad dB'\in \Delta\iff k\mid d.\end{equation}  The proofs of the two statements are the same, so we focus on the first. In one direction, it is easy to check that
\[nA'=A=p_{*}(\tilde{A})\in \Delta.\]  For the converse, suppose that $dA' = p_*(C)$ for some class $C\in  \Num(\tilde{X})$. Then, since $p^*(A)=n\tilde{A}$, and $A\cdot A=0$, it follows that $\tilde{A}\cdot C=0$. A result of  Kani\footnote{Proof of Proposition 2.3 in E. Kani, Elliptic curves on Abelian surfaces, Manuscripta Math. 84, 199 -- 223 (1994).} then shows that $C$ is a multiple of $\tilde{A}$, which proves the claim.

\smallskip

 Suppose now that we have a FM transform $\Phi\colon D(Y)\to D(X)$. We must show that $Y\isom X$. By the existence of the lift of $\Phi$ to the canonical covers, and Cor. \ref{abelsur}, we know that  $\Phi$ is a sheaf transform. Composing with a shift we can therefore suppose that  $\Phi$ takes skyscrapers of points on $Y$ to sheaves on $X$ of class \[v=[\Phi(\OO_y)]=(r,aA'+bB',s)\in \Z\oplus\Num(X)\oplus\Z.\]  We have  $r\geq 0$. The condition $\chi(v,v)=0$ gives  $rs=ab$. By the existence of the lift, $\Phi(\OO_y)$ is the push-forward of an sheaf on $\tilde{X}$. It follows that $n\mid  r$ and $aA'+bB'\in \Delta$.
 
 Suppose first that  $r=0$. Then  either $a=0$ or $b=0$. It follows that $\Phi(\OO_y)$ is supported on a fibre of one of the two elliptic fibrations $\pi_A,\pi_B$. Let us denote this fibration by $\pi\colon X\to C$. 
 The argument of Prop. \ref{elika} then applies, and shows that $Y=J_{X/C}(b)$, with $b$ coprime to $\lambda=\lambda_{X/C}$. 
 Since $A'\cdot B=k$ and $A\cdot B'=n$, we must have \[\lambda\in \{1,2,3,4,6\}.\] But  the only units in $\Z/\lambda \Z$ are then $\pm 1$,  so it follows from the arguments of Remark \ref{misc} that $Y\isom X$.
  
 Suppose then that $r>0$. Consider the elliptic fibration \[\pi_A\colon X\to A/J\] whose fibres have  class $B$. The divisor $A'$ gives a $k$-section. Sheaves of class $v$ have  rank $r$ and fibre degree $ka$. Take $x,y\in \Z$ with \[xr +yka= h: =\operatorname{gcd}(r,ka).\]
  We claim that we can suppose that $k\mid x$.
 Assuming this, consider
\[\mat{ka/h}{-r/h}{x}{y} \col{r}{ka}=\col{0}{h}.\]
The negative of the matrix on the left is then as in Theorem \ref{super}, so  lifts to a relative Fourier-Mukai transform. Since all relative Jacobians of $X$ are isomorphic to $X$ as above, we obtain an equivalence  \[\Psi\colon D(X)\to D(X)\] which takes $v$ to a vector of rank $0$. Applying the argument above to $\Psi\circ \Phi$ we conclude that $Y\isom X$, which completes the proof.

To prove the claim we can  assume that  $k\in\{2,3\}$ is  prime. Writing $r=hr'$ and $ka=ht'$ we are looking for pairs $x,y\in \Z$ with $xr'+yt'=1$. If $(x,y)$ is such a pair, so is $(x+t',y-r')$. Thus we can find such a pair with $k\mid x$ precisely if $k\nmid t'$. So what  we must prove is that the highest power of $k$ which divides $ka$  also divides $r$. 
 
 Consulting the list of possibilities, we either have $n=k\in\{2,3\}$ or $(n,k)=(4,2)$.   Recall that \[aA'+bB'\in \Delta, \qquad rs=ab, \qquad n\mid r.\]  Clearly, one of $a$ or $b$ is divisible by $k$. In the cases when $n=k$, using  \eqref{err_stuff}, it follows  that since one of the two classes $aA',bB'$ lie in $\Delta$, they both do, and hence $k$ divides both $a$ and $b$. In the remaining case $(n,k)=(4,2)$ we have $4|ab$, and we again conclude that $k$ divides both $a$ and $b$. The primitivity of $v$  now ensures that $k\nmid s$, and the claim follows.

%
%\small
%
%Department of Mathematics and Statistics, The University of Edinburgh,
%King's Buildings, Mayfield Road, Edinburgh, EH9 3JZ, UK.
%
%\smallskip
%email: {\tt tab@maths.ed.ac.uk}
%{\tt\ \ \ \ \ \ \ \ \ \ ama@maths.ed.ac.uk}

\end{document}